\documentclass[12pt]{article}
\usepackage[mathscr]{euscript}
\usepackage{amsmath,amssymb}
\usepackage{verbatim}
\usepackage{dsfont,bm}
\usepackage{color}
\usepackage{graphicx}
\usepackage{amsthm}
\usepackage{enumerate}
\usepackage{esint} %\fint

\usepackage{xcolor}
\definecolor{darkgreen}{rgb}{0,0.75,0}
\definecolor{darkred}{rgb}{0.75,0,0}
\definecolor{darkmagenta}{rgb}{0.5,0,0.5}
\usepackage[colorlinks,citecolor=darkgreen,linkcolor=darkred,urlcolor=darkmagenta]{hyperref}
\newtheorem{theorem}{Theorem}[section]

\newtheorem{lemma}[theorem]{Lemma}

\newtheorem{prop}[theorem]{Proposition}
\newtheorem*{notation}{Notation}

\theoremstyle{definition}
\newtheorem{definition}[theorem]{Definition}

\newtheorem{remark}[theorem]{Remark}

\numberwithin{equation}{section}

\def\be{\begin{equation}}
\def\ee{\end{equation}}
\def\bes{\begin{equation*}}
\def\ees{\end{equation*}}

\newcommand{\mr}[1]{{\tt \href{http://www.ams.org/mathscinet-getitem?mr=#1}{MR#1}}}

\newcommand{\set}[1]{\left\{ #1 \right\}}

\newcommand{\abs}[1]{{\left\vert\kern-0.25ex #1
		\kern-0.25ex\right\vert}}
\newcommand{\one}{\mathds{1}} %indicator

\newcommand{\diag}[0]{\operatorname{diag}}
\DeclareMathOperator*{\esssup}{ess\,sup}
\DeclareMathOperator*{\essinf}{ess\,inf}

 \def\sE {{\mathcal E}} \def\sF {{\mathcal F}}

\def\sM {{\mathcal M}} \def\sN {{\mathcal N}}

 \def\bE {{\mathbb E}}

\def\bP {{\mathbb P}}  \def\bR {{\mathbb R}}

\def\sms{\smallskip}

\def\ignore#1{}

 \def\Lam {\Lambda} 
\def\eps{\varepsilon}

\def\to {\rightarrow}

\def\q{\quad} \def\qq{\qquad}
\def\dint{\int\kern-.6em\int}

\newcommand\restr[2]{{% we make the whole thing an ordinary symbol
		\left.\kern-\nulldelimiterspace % automatically resize the bar with \right
		#1 % the function
		\vphantom{\big|} % pretend it's a little taller at normal size
		\right|_{#2} % this is the delimiter
	}} %restriction of a function
	\def\wt{\widetilde}
	\def\wh{\widehat}
	
	\def\be{\begin{equation}}
	\def\ee{\end{equation}}
	\def\bes{\begin{equation*}}
	\def\ees{\end{equation*}}
	\def\ba{\begin{align}}
	\def\ea{\end{align}}
	\def\xxea{\end{align}}
\def\bas{\begin{align*}}
\def\eas{\end{align*}}

\def\proof{{\smallskip\noindent {\em Proof. }}}
\def\qed{{\hfill $\square$ \bigskip}}

\definecolor{dgreen}{rgb}{0, 0.6, 0.1}
\definecolor{dblue}{rgb}{0, 0.0, 0.6}
\definecolor{vdblue}{rgb}{0,.08, 0.45}
\definecolor{dred}{rgb}{0.7, 0.0, 0.0}
\definecolor{vdblue}{rgb}{0,.08, 0.45}

\definecolor{purple}{rgb}{0.6, 0.0, 0.6}
\definecolor{mytext}{rgb}{0.1, 0.1, 0.1}

\begin{document}
	
	\font\titlefont=cmbx14 scaled\magstep1
	\title{Parabolic Harnack inequality implies the existence of jump kernel \titlefont }    
	\author{	Guanhua Liu\footnote{Research partially supported by China Scholarship Council.} ,	Mathav Murugan\footnote{Research partially supported by NSERC (Canada).} 
		}
	\maketitle
	\vspace{-0.5cm}
%\ver
\begin{abstract}
We prove that the parabolic Harnack inequality implies the existence of jump kernel for symmetric pure jump process. This allows us to remove a technical assumption on the jumping measure in the recent characterization of the parabolic Harnack inequality for pure jump processes by Chen, Kumagai and Wang. The key ingredients of our proof are the L\'evy system formula and a near-diagonal heat kernel lower bound.
\vskip.2cm
%\noindent {\it Keywords:} 
\
\end{abstract}

\section{Introduction}
The parabolic Harnack inequality is a fundamental regularity estimate for non-negative solutions to the heat equation and its variants. Important applications of the parabolic Harnack inequality are apriori H\"older regularity of solutions, the existence of heat kernel, and  bounds on the heat kernel. We refer to the survey \cite{Kas} for an introduction to Harnack inequalities and variants.

A major result on the parabolic Harnack inequality is its characterization by simpler geometric and analytic properties. In the context of diffusions on Riemannian manifolds this characterization was established by Grigor'yan \cite{Gri} and Saloff-Coste \cite{Sal} for the classical space-time scaling (time scales like square of space). This characterization was extended and modified to many settings including diffusions on metric measure spaces \cite{Stu}, nearest neighbor walks on graphs \cite{Del}, and for anomalous space time scaling \cite{BB,BBK06,GHL} by several authors. A similar characterization of the parabolic Harnack inequality for jump processes remained open until a recent breakthrough by Chen, Kumagai and Wang \cite{CKW2}. 

The purpose of this note is to show that the parabolic Harnack inequality implies the existence of the jump kernel for pure jump processes. In other words, we show that the jumping measure is absolutely continuous with respect to the product measure $\mu \otimes \mu$, where $\mu$ is the symmetric (reference) measure for the jump process. As a consequence, we remove a technical hypothesis on the jumping measure assumed in \cite{CKW2} for characterizing the parabolic Harnack inequality 
(see Remark \ref{r:ckw}).

\subsection{Framework and result}
Let $(M,d)$ be a complete, locally compact, separable metric space, and let $\mu$ be a positive Radon measure on $M$ with full support. Such a triple $(M,d,\mu)$ is called a \emph{metric measure space}.
 We assume that $(M,d)$ is \emph{unbounded}; that is, $\sup_{x,y \in M}d(x,y)=\infty$.
  We set $B(x,r):=\{y\in X\mid d(x,y)<r\}$ and $V(x,r)= \mu(B(x,r))$ for $x \in M, r \in (0,\infty)$. We assume that 
the measure $\mu$ satisfies the following \emph{volume doubling} property \eqref{VD}: there exists $C_D >1$ such that
\begin{equation}\label{VD} \tag{VD}
V(x,2r)\le C_{D} V(x,r), \q \mbox{for any $x \in M, r \in (0,\infty)$}.
\end{equation}

We consider a \emph{symmetric Dirichlet form} $(\sE,
\sF)$ on $L^2(M,\mu)$. In other words, $\sF$ is a dense linear subspace of $L^2(M,\mu)$, $\sE:\sF \times \sF \to \bR$ is symmetric, non-negative definite, bilinear form that is \emph{closed} ($\sF$ is a Hilbert space under the inner product $\sE_1(\cdot,\cdot)= \sE(\cdot,\cdot) + \langle \cdot,\cdot \rangle_{L^2(M,\mu)}$) and \emph{Markovian} (for any $f \in \sF$, we have $\wh f:=(0\vee f)\wedge 1 \in \sF$ and $\sE(\wh f,\wh f) \le \sE(f,f)$). We assume that $(\sE,\sF)$ is \emph{regular}; that is, $\mathcal{F}\cap\mathcal{C}_{\mathrm{c}}(X)$ is dense both in $(\mathcal{F},\mathcal{E}_{1})$
and in $(\mathcal{C}_{\mathrm{c}}(X),\|\cdot\|_{\mathrm{sup}})$.  We assume that $(\sE,\sF)$ is a \emph{pure jump} type Dirichlet form; that is, there exists a symmetric positive Radon measure on $M \times M \setminus \diag$ such that 
\[
\sE(f,f) = \int_{M \times M \setminus \diag} (f(x)-f(y))^2\, J(dx,dy), \q \mbox{for all $f \in \sF$,}
\]
where $\diag = \set{(x,x) \mid x \in M}$ denotes the diagonal. The Radon measure $J$ is called the \emph{jumping measure}; cf.~\cite[Theorem 3.2.1]{FOT}. 
We say that the Dirichlet form $(\sE,\sF)$ on $L^2(X,\mu)$ \emph{admits a jump kernel} if $J$ is absolutely continuous with respect to the product measure $\mu \otimes \mu$ on $M \times M \setminus \diag$. If the Dirichlet form admits a jump kernel, then the Radon-Nikodym derivative of $J$ with respect to $\mu \otimes \mu$ is called the \emph{jump kernel}. In other words, (if it exists) the jump kernel $j(\cdot,\cdot)$  is a measurable function such that $J(dx,dy)=j(x,y) \mu(dx) \mu(dy)$. The  central question of this work whether or not a pure jump type, regular Dirichlet form admits a jump kernel.

Every regular Dirichlet form  $(\sE,
\sF)$ on $L^2(M,\mu)$ has an associated $\mu$-symmetric \emph{Hunt process}
$X= \set{X_t, t \ge0, \bP^x, x\in M \setminus \sN}$, where $\sN$ is a \emph{properly exceptional} set for $(\sE,\sF)$; that is, $\mu(\sN)=0$ and $\mathbb{P}^x(X_t\in\mathcal{N} \mbox{ for some $t>0$})=0$.
This Hunt process is unique up to the choice of a properly exceptional set \cite[Theorems 4.2.8 and 7.2.1]{FOT}. 
% We denote the expectation with respect to $\bP^x$ by $\bE^x$.
	Let $Z_t=(V_t,X_t)$ be the associated space-time process ($\mathbb R \times M$-valued process) defined by $V_t=V_0-t$. The law of the space time process $s \mapsto Z_s$ starting from $(t,x)$ will be denoted by $\bP^{(t,x)}$. The  expectation with respect to $\bP^{(t,x)}$ is denoted by $\bE^{(t,x)}$.
	We say that $A \subset [0,\infty) \times M$ is \emph{nearly Borel measurable} if for any Borel probability measure $\mu_0$ on $[0,\infty) \times M$, there are Borel measurable subsets $A_1,A_2$ such that $A_1 \subset A \subset A_2$ and satisfies $\bP^{\mu_0}(Z_t \in A_2 \setminus A_1 \mbox{ for some $t\ge 0$})=0$. The collection of nearly Borel measurable subsets of $[0,\infty) \times M$ forms a $\sigma$-field, which is called nearly Borel measurable $\sigma$-field.
	We recall the (probabilistic) definition of the parabolic Harnack inequality.
\begin{definition} {\rm

	We say that a nearly Borel measurable function $u:[0,\infty) \times M \to \mathbb{R}$ is \emph{caloric} on $D= (a,b) \times B(x_0,r)$ for the Markov process $X$ if there is a property exceptional set $\mathcal{N}_u$ of the Markov process $X$ such that for any relatively compact open subset $U$ of $D$, we have
	\[
	u(t,x)= \mathbb{E}^{(t,x)} u (Z_{\tau_U}), \quad \mbox{for all $(t,x) \in U \cap \left( [0,\infty) \times (M\setminus \mathcal{N}_u)\right)$.}
	\]
	
	Let $\phi:[0,\infty) \to [0,\infty)$ be a homeomorphism (and hence strictly increasing with $\phi(0)=0$). We say that the parabolic Harnack inequality \ref{PHI} holds for the process $X$, if  there exist constants $c_0 \in (0,1)$, $0<C_1<C_2<C_3<C_4$ and  $C_5>1$ such that for any $x_0\in M, t_0\ge 0, r>0$, provided that $u:\mathbb{R}_+\times M\to\mathbb{R}_+$ is caloric on $(t_0,t_0+C_4\phi(r))\times B(x_0,r)$, we always have
	\begin{equation}\label{PHI} \tag{$\operatorname{PHI(\phi)}$}
	\esssup\limits_{(t_0+C_1\phi(r),t_0+C_2\phi(r))\times B(x_0,c_0r)}u\le C_5 \essinf\limits_{(t_0+C_3\phi(r),t_0+C_4\phi(r))\times B(x_0,c_0r)}u.
	\end{equation}
}
\end{definition}
The main result of our work is that the parabolic Harnack inequality implies the existence of jump kernel.
\begin{theorem} \label{t:jkexists}
Let $(M,d,\mu)$ be an  unbounded, complete, separable, locally compact metric measure space, where $\mu$ is a Radon measure with full support on $(M,d)$ that satisfies the volume doubling property \eqref{VD}. Let $(\sE,\sF)$ be a symmetric Dirichlet form on $L^2(M,\mu)$ of pure jump type and let $X$ be the corresponding $\mu$-symmetric Hunt process. Let $\phi:[0,\infty) \to [0,\infty)$ be a homeomorphism such that there exist constants $C_\phi \ge 1$, $\beta_2 \ge \beta_1 >0$ such that 
\be  \label{e:regphi}
C_\phi^{-1} \left(\frac R r\right)^{\beta_1} \le \frac{\phi(R)}{\phi(r)} \le C_\phi \left(\frac R r\right)^{\beta_2} \q \mbox{for all $0 < r \le R$.}
\ee
If the process $X$ satisfies the parabolic Harnack inequality \ref{PHI}, then the Dirichlet form $(\sE,\sF)$ on $L^2(M,\mu)$ admits a jump kernel. 
\end{theorem}
\begin{remark}\label{r:ckw} {\rm 
		\begin{enumerate} 
\item[(a)]		Let $J(dx,dy)$ denote the jumping measure for $(\sE,
		\sF)$. Assume that there is a kernel $\wt J(x,dy)$ (in other words, $x \mapsto \wt J(x,A)$ is a Borel measurable function for any Borel set $A$, and that $A \mapsto \wt J(x,A)$ is a Borel measure on $M$ for any $x \in M$) such that
		\be \label{e:tech}
		J(dx, dy) = \wt J(x,dy) \mu(dy). 
		\ee
 Theorem \ref{t:jkexists} was shown under the additional assumption that a kernel $\wt J(x,dy)$ exists and satisfies \eqref{e:tech}  in \cite[Proposition 3.3]{CKW2} (see also \cite[Proposition 4.7]{BBK09} for a similar result and proof).  This assumption can be viewed as a weaker form of the existence of jump kernel and was assumed throughout \cite{CKW2}.
As a consequence of Theorem \ref{t:jkexists}, we could remove the assumption \eqref{e:tech} in the characterization of parabolic Harnack inequality in \cite{CKW2}. 
\item[(b)] As explained in \cite[Remark 1.22]{CKW2}, the condition that the metric space is unbounded can be relaxed. Our proof of Theorem \ref{t:jkexists} also extends to the case where there are non-zero diffusion and jump parts as considered in \cite{CKW3}. 
\end{enumerate}
}
\end{remark}
In the proof of Theorem \ref{t:jkexists}, we consider the same caloric function used in \cite[Proposition 3.3]{CKW2} and   \cite[Proposition 4.7]{BBK09}. However, the argument in \cite{CKW2} requires a L\'evy system formula (see \cite[Lemma 2.11]{CKW2}) that relies on the assumption \eqref{e:tech}. To overcome the difficulty, we use a more abstract L\'evy system formula that does not rely on \eqref{e:tech}.
The main new ingredient in our proof is the use of a near diagonal lower bound on the heat kernel to obtain useful quantitative estimates on the jumping measure. In particular, we use \emph{both} upper and lower bounds on the heat kernel while the argument in \cite[Proposition 3.3]{CKW2} uses only upper bound on the heat kernel.

\begin{notation}
	Throughout this paper, we use the following notations and conventions.
\begin{enumerate}[(a)]
	\item For a measurable function $f \ge 0$ and a measure $\mu$, by $f \cdot \mu$, we denote the measure $A \mapsto \int_A f\,d\mu$.
	\item For a measure $\mu$ and a function $f$, the integral $\int f \,d \mu$ is denoted by $\langle \mu , f \rangle$.
	\item The notation  $A \lesssim B$ for quantities $A$ and $B$ indicates the existence of an
	implicit constant $C \ge 1$ depending on some inessential parameters such that $A \le CB$. We write $A \asymp B$, if $A \lesssim B$ and $B \lesssim A$.
%	\item Every function in the domain of a regular Dirichlet form is  represented by a quasi-continuous function.
%	\item For a topological space $M$, by $\sB(M)$ we denote the $\sigma$-field of Borel sets.
\end{enumerate}

\end{notation}

\section{Proof}
The proof of Theorem \ref{t:jkexists} relies on two key ingredients. The ingredients are  bounds on the heat kernel and a L\'evy system formula, which we recall in \textsection  \ref{s:nle} and \textsection \ref{s:levysystem} respectively.
After these preliminaries, we present the proof of Theorem \ref{t:jkexists} in \textsection \ref{s:jkexists}.

\subsection{Heat kernel} \label{s:nle}
%***Define heat kernel, PHI implies existence of heat kernel and NLE (for Dirichlet heat kernel).

We recall the notion of \emph{heat kernel}.
Let $(M,d,\mu)$ be a metric measure space and let  $(\mathcal{E},\mathcal{F})$ be a regular Dirichlet form on $L^2(M,\mu)$. Let $(X_t, t \ge 0, \mathbb{P}_x, x \in  M \setminus \mathcal{N})$ be the corresponding $\mu$-symmetric Hunt process, where $\mathcal N$ is a properly exceptional set for $(\sE,\sF)$.  Let $\set{P_t}$ note the corresponding \emph{Markov semigroup} \cite[Theorem 1.4.1]{FOT}.
The \emph{heat kernel} associated with the Markov semigroup $\{P_t\}$ (if it exists) is a family of measurable functions $p(t,\cdot,\cdot):M  \times M  \mapsto [0,\infty)$ for every $t > 0$, such that
\begin{align}
P_tf(x) &= \int p(t,x,y)  f(y)\, \mu(dy), \q \mbox{for all $f \in L^2(M,\mu), t>0$ and $x \in M$,} \label{e:hkint} \\
p(t,x,y) &= p(t,y,x) \q \mbox{for all $x, y \in M$ and $t>0$,} \label{e:hksymm} \\
p(t+s,x,y)&= \int p(s,x,y) p(t,y,z) \,\mu(dy), \q \mbox{for all $t,s>0$ and $x,y \in M$.} \label{e:ck}
\end{align}
For an open set $B$, let $X^B$ denote the $\mu$-symmetric Hunt process on $B$ obtained from $X$ killed upon exiting $B$ \cite[Theorems 3.3.8 and 3.3.9]{CF}; that is, 
\[ X^B(t)= \begin{cases}
X(t) & \mbox{for $t < \tau_B$,} \\
\Delta & \mbox{for $t \ge \tau_B$}
\end{cases},\]
where $\Delta$ denotes the cemetery state and $\tau_B= \inf\set{s>0 \mid X_t \notin B}$ denote the exit time of $B$. Let $\set{P_t^B}$ denote the Markov semigroup on $L^2(B,\mu)$.
   The \emph{heat kernel} associated with the Markov semigroup $\{P_t^B\}$ (if it exists) is denoted by $p^B(t,\cdot,\cdot)$. We recall the existence and  bounds on the heat kernels for $X$ and $X^B$ from \cite{CKW2}.

\begin{prop} \label{p:hke} \cite[Propositions 3.1 and 3.2]{CKW2}
	Let $(M,d,\mu)$ be an  unbounded, complete, separable, locally compact metric measure space, where $\mu$ is a Radon measure with full support on $(M,d)$ that satisfies the volume doubling property \eqref{VD}. Let $(\sE,\sF)$ be a symmetric Dirichlet form on $L^2(M,\mu)$ of pure jump type and let $X$ be the corresponding $\mu$-symmetric Hunt process. Let $\phi:[0,\infty) \to [0,\infty)$ be a homeomorphism such that there exist constants $C_\phi \ge 1$, $\beta_2 \ge \beta_1 >0$ satisfying \eqref{e:regphi}.
	Assume further that  $X$ satisfies the parabolic Harnack inequality \ref{PHI}. Then
	\begin{enumerate}[(a)]
		\item The process $X$ has a continuous heat kernel $p:(0,\infty) \times M \times M \to [0,\infty)$ that satisfies the following upper bound. There exists a constant $C_{U}>0$ and a properly exceptional set $\sN$ for $X$ such that,
		\begin{equation}\label{UE}
		p_t(x,y)\le\frac{C_{U}}{V(x,\phi^{-1}(t))}, \q \mbox{for all $x,y\in M \setminus \sN$ and for all $t>0$,} \tag{$\operatorname{UHKD(\phi)}$}
		\end{equation}
		
		\item For every ball $B=B(x_0,r)$, let $X^B$ denote the process obtained from $X$ killed upon exiting $B$. Then $X^B$ has a heat kernel $p^B:(0,\infty) \times B \times B \to [0,\infty)$ and satisfies the following lower bound:
 there exists $c_L>0, \delta_N \in (0,1)$ and a properly exceptional set $\sN$ for $X$ such that for any $x_0 \in M, r<0, 0<t \le \phi(\delta_N r)$ and $B=B(x_0,r)$,
		\begin{equation}\label{NLE} \tag{$\operatorname{NDL(\phi)}$}
		p_t^{B(x_0,r)}(x,y)\ge\frac{c_{L}}{V(x_0,\phi^{-1}(t))},\quad \mbox{  for all $x,y\in B(x_0,\delta_N\phi^{-1}(t)) \setminus \sN$}.
		\end{equation}
		
	\end{enumerate}
\end{prop}
%IS it possible to remove $\sN$ using parabolic Holder regularity?
%REMARK: properly exceptional set has capacity zero, so is not charged by any Revuz measure.
\begin{remark} \label{r:conthk}{\rm 
		\begin{enumerate}[(a)]
			\item We remark that the proofs of  \cite[Propositions 3.1 and 3.2]{CKW2} do not rely on the assumption \eqref{e:tech} or the reverse volume doubling property.
			\item Using a standard parabolic H\"older regularity estimate (see \cite[Corollary 4.5 and Lemma 4.6]{BGK}), we may assume that $(t,x,y) \mapsto p_t(x,y)$ is continuous on $(0,\infty) \times M \times M$. Similarly, for any ball $B(x_0,r)$, we may assume that $(t,x,y) \mapsto p_t^{B(x_0,r)}(x,y)$ is continuous in $(0,\infty)\times B(x_0,r) \times B(x_0,r)$. In particular, we may assume that the exceptional set $\sN$ in the estimates \ref{UE} and \ref{NLE} is the empty set.
		\end{enumerate}	
}
\end{remark}

\subsection{L\'evy system formula} \label{s:levysystem}
In this section, we collect some useful facts on the L\'evy system formula and positive continuous additive functionals.

Consider a  $\mu$-symmetric \emph{Hunt process}
$X= \set{\Omega, \sM, X_t, t \ge0, \bP^x}$, where $\sN$ is a  properly exceptional set for $(\sE,\sF)$ on $L^2(X,\mu)$ and $(\Omega, \sM,\bP^x)$. For any measure $\nu$ on $M$, we denote by $\bP^\nu$ the measure $\bP^\nu (A)= \int_{M} \bP^x (A)\, d\nu(x)$. Any function $f$ on $M$ is extended to $M_\partial:= M \cup \set{\Delta}$ by setting $f(\Delta)=0$, where $\Delta$ denotes the cemetery state. The set $M_\partial$ as a topological space is the one point compactification of $M$. Let $(\sM_t)_{0 \le t\le \infty}$ denote the minimum augmented admissible filtration on $\Omega$.

A collection of random variables $A:= \{A_s:\Omega\to\mathbb{R}_+|s\in\mathbb{R}_+\}$, is called a \emph{positive continuous additive functional} (for short, a PCAF), if it satisfies the following conditions:
\begin{enumerate}[(i)]
	\item  $A_t(\cdot)$ is $(\sM_t)$-measurable,
	\item there exist a set $\Lam \in \sM_\infty$ and  an exceptional set $\sN \subset M$ for $X$ such that $\bP^x(\Lambda)=1$
 for all $x \in M \setminus \sN$ and $\theta_t \Lam \subset \Lam$ for all $t>0$, where $\theta_t$ denotes the shift map on $\Omega$.
	\item  For any $\omega\in\Lambda$, $t\mapsto A_t(\omega)$ is continuous, non-negative with $A_0(\omega)=0$, $A_t(\omega)=A_{\zeta(\omega)}(\omega)$ for $t\ge\zeta(\omega)$, and $A_{t+s}(\omega)=A_t(\omega)+A_s(\theta_t\omega)$ for any $s,t\ge 0$. Here  $\zeta(\cdot)$ denotes the life time of the process.
\end{enumerate}
The sets $\Lam$ and $\sN$ are referred to as a \emph{defining set} and \emph{exceptional set} of the PCAF $A_t$ respectively. If $\sN$ can be taken to the empty set, then we say that $A_t$ is a PCAF in the \emph{strict sense}.

% Given a PCAF $A$ and a Borel function $f\ge 0$, we denote $(fA)_t:=\int_0^tf(X_s)\,dA_s$, which is still a random variable.

 A measure $\nu$ is called the \emph{Revuz measure} of the PCAF $A$, if and only if for any non-negative Borel functions $h$ and $f$,
\begin{equation}\label{Revuz}
\mathbb{E}^{h\cdot\mu} \left(\int_0^t f(X_s(\omega))\,dA_s(\omega) \right)=\int_0^t\langle f\cdot\nu,P_s h\rangle \, ds,
\end{equation}
where $P_s$ denotes the Markov semigroup corresponding to the Hunt process.
By \cite[Theorem 5.1.4]{FOT}, the Revuz measure $\nu$ is uniquely determined by $A$ and does not charge any set of zero capacity. In particular,
\be \label{e:smoothnocharge}
\nu(\sN)=0, \q \mbox{for any properly exceptional set $\sN$.}
\ee

Every Hunt process has a \emph{L\'evy system} $(N,H)$ \cite[Appendix A.3.4]{CF}. Recall that 
a pair $(N,H)$ is a L\'evy sytem for the Hunt process $X$ if $N(x,dy)$ is a kernel on $M_\partial$ equipped with the Borel $\sigma$-field   and $H$ is a PCAF in the strict sense  satisfying the following property: for any non-negative Borel function $F: M_\partial \times M_\partial \to [0,\infty)$ such that $F(x,x)=0$ for all $x \in M_\partial$, we have
\be \label{e:LSbasic}
 \mathbb{E}^z\left[\sum\limits_{s\le t}F(X_{s-},X_s)\right]=\mathbb{E}^z\left[\int_0^t\int_MF(X_s,y)N(X_s,dy)dH_s\right].
 \ee
  The property in \eqref{e:LSbasic} called the \emph{L\'evy system formula} and admits the following generalization. By \cite[(A.3.33)]{CF}, for any non-negative Borel function $g$ on $(0,\infty)$, any $z \in M_\partial$,  any $(\sM_t)$-stopping time $T$, and any non-negative Borel function $F: M_\partial \times M_\partial \to [0,\infty)$ such that $F(x,x)=0$ for all $x \in M_\partial$, we have
 \be \label{e:LS}
 \mathbb{E}^z\left[\sum\limits_{0<s\le  T}g(s)F(X_{s-},X_s)\right]=\mathbb{E}^z\left[\int_0^ T g(s) \int_M F(X_s,y)N(X_s,dy)dH_s\right].
 \ee
By \cite[(5.3.6) and Theorem 5.3.1]{FOT}, we know if $\nu$ is the Revuz measure of $H$, then 
\begin{equation} \label{e:levy-jump}
J(dx,dy)= \frac 1 2 N(x,dy)\nu(dx).
\end{equation}

\begin{lemma}\cite[Proposition 4.1.10]{CF}
	Let $H$ be a PCAF for the process $(X_t)$ and let $\nu$ be the corresponding Revuz measure. For any open set $D$ the process $(H_{t\wedge \tau_D})$ is a PCAF for the  process $X^D$ killed upon exiting $D$ and its Revuz measure is $\nu_D(\cdot)=\nu(D\cap \cdot)$, where $\tau_D=\inf\set{t>0: X_t \notin D}$. In particular, we have
	\begin{equation}\label{Revuz-loc}
	\mathbb{E}^{h\cdot\mu}\left( \int_0^{\tau_D\wedge t} f(X_s)\,dH_s \right)=\int_0^t\langle f\cdot\nu_D,P_s^{D}h\rangle \, ds,
	\end{equation}
	for all non-negative measurable functions $f,h:D \to [0,\infty)$, where $P^D_s$ denotes the Markov semigroup corresponding to the $X^D$.
\end{lemma}

\subsection{Existence of jump kernel} \label{s:jkexists}

%We recall a parabolic H\"older regularity estimate for caloric functions. This is a well known consequence of the parabolic Harnack inequality. We refer to \cite[Proposition 4.4 and Corollary 4.5]{BGK} for a proof.
%\begin{lem} \label{l:phr} 	Let $(M,d,\mu)$  metric measure space equipped with    a symmetric Dirichlet form  $(\sE,\sF)$ on $L^2(M,\mu)$ of pure jump type satsifying the assumptions of Proposition \ref{p:hke}.  Let $\phi:[0,\infty) \to [0,\infty)$ be a homeomorphism such that there exist constants $C_\phi \ge 1$, $\beta_2 \ge \beta_1 >0$ satsifying \eqref{e:regphi}.
%	Assume further that the Hunt process $X$ corresponding to the above Dirichlet form satisfies the parabolic Harnack inequality \ref{PHI}. 
%Let $u:[0,\infty) \times M$ be caloric in $[t_0, t_0 + \phi(r)] \times B(x_0,r)$, for some $t_0 \ge 0$.
%\end{lem}

The following estimate on the jump kernel plays a crucial role in the proof of Theorem \ref{t:jkexists}. This estimate can be viewed as an integrated version of the condition (UJS) considered in \cite[Definition 1.18]{CKW2}.
\begin{lemma} \label{l:jumpav}
 Under the assumptions of Theorem \ref{t:jkexists}, there exists $\delta\in (0,1), C_J>0$ such that for any pair of balls $B_i=B(x_i,r_i), i=1,2$ with $d(x_1,x_2)>r_1+r_2$ and for any ball $B'= B(x',r') \subset B(x_1,\delta r_1)$ such that  $r' \le \delta r_1$, we have the estimate
\[
 J(B'\times B_2) \le C_J \frac{\mu(B')}{\mu(B_1)} J(B_1 \times B_2).
\]
\end{lemma}
\proof 
Let $c_0 \in (0,1), C_1,C_2,C_3,C_4 >0, C_5 >1$ denote the constants in \ref{PHI}.
Let $B_i=B(x_i,r_i), i=1,2$ be balls such that $r_1+r_2 < d(x_1,x_2)$.
Set $f_h(t,z)= \one_{(C\phi(r)-h, C\phi(r))}(t) \one_{B_2}(z)$, where $C=(C_1+C_2)/2$ and $h \in (0, C\phi(r))$. Then
\begin{equation}
u_h(t,x)= \begin{cases}
\mathbb{E}^x[f_h(t-\tau_{B_1},X_{\tau_{B_1}}); \tau_{B_1} \le t] & \mbox{if $x \in M \setminus \mathcal{N}$, $t>0$}\\
0 & \mbox{if $x \in \mathcal{N}$, $t>0$,}
\end{cases}
\end{equation}
is caloric in $(0,\infty) \times B_1$, where $\sN$ is an exceptional set for the corresponding Hunt process $X$ and $\tau_{B_1} = \inf \set{t>0 \mid X_t \notin B_1}$ denote the exit time from $B_1$. By Remark \ref{r:conthk}, we may assume that the heat kernel corresponding to the process killed upon exiting $B_1$ given by $(t,x,y) \mapsto p^{B_1}_t(x,y)$ is continuous in $(0,\infty) \times B_1 \times B_1$.

We choose a L\'evy system $(N,H)$ for the process $X$. 
Let  $\nu$ denote the Revuz measure of $H$, where $H$ is a PCAF in the strict sense.
Set  $g(x)=N(x,B_2)$.

%Let $G_{t,h}: (0,\infty) \times M \times M \to [0,\infty)$ be defined as
%\[
%G_{t,h}(s,x,y)=  \one_{(t-C\phi(r),t-C\phi(r)+h)}(s) \one_{B_1}(x) \one_{B_2}(y).
%\]
For any $t > C \phi(r)$, for quasi-every $x \in B_1$, for any $h \in (0,C\phi(r))$, and for any $s_1 \in (0,t-C\phi(r))$), we have
\begin{align}
u_h(t,x) &= \bE^x \left[\sum_{s\le \tau_{B_1}} \one_{(t-C\phi(r),t-C\phi(r)+h)}(s) \one_{B_1}(x) \one_{B_2}(y) \right] \nonumber \\
&= \bE^x \left[ \int_0^{\tau_{B_1}} \int_M \one_{(t-C\phi(r),t-C\phi(r)+h)}(s) \one_{B_1}(x) \one_{B_2}(y)\, N(X_s,dy)\, dH_s \right] \q \mbox{(by \eqref{e:LS})} \nonumber \\
&= \bE^x \left[ \int_{(t-C\phi(r)) \wedge \tau_{B_1}}^{((t- C\phi(r)+h) \wedge \tau_{B_1})} N(X_s,B_2)\, dH_s\right] = \bE_x \left[ \int_{(t-C\phi(r)) \wedge \tau_{B_1}}^{((t- C\phi(r)+h) \wedge \tau_{B_1})} g(X_s)\, dH_s\right]\nonumber \\
&= \bE^{p_{s_1}^{B_1}(x,\cdot)\cdot\mu} \left[ \int_{(t-C\phi(r)- s_1) \wedge \tau_{B_1}}^{((t- C\phi(r)+h-s_1) \wedge \tau_{B_1})} g(X_s)\, dH_s\right] \q \mbox{(by Markov property)}\nonumber\\
&= \int_{t- C \phi(r)-s_1}^{t- C \phi(r)+h-s_1} \langle g\cdot \nu_{B_1}, P_s^{B_1}p_{s_1}^{B_1}(x,\cdot) \rangle \, ds \q \mbox{(by \eqref{Revuz-loc})} \nonumber\\
&= \int_{t- C \phi(r)-s_1}^{t- C \phi(r)+h-s_1} \langle g\cdot \nu_{B_1}, p_{s+s_1}^{B_1}(x,\cdot) \rangle \, ds\q \mbox{(by \eqref{e:hkint} and \eqref{e:ck})} \nonumber\\
%&= \int_{t- C \phi(r)}^{t-C\phi(r)+h} \int_{B_1} p_s^{B_1}(x,w) N(w,B_2)\, \nu_{B_1}(dw)\, ds \nonumber\\
&= \int_{t- C \phi(r)}^{t-C\phi(r)+h} \int_{B_1} p_s^{B_1}(x,w) N(w,B_2)\, \nu(dw)\, ds \q \mbox{(since $g(\cdot)=N(\cdot,B_2)$)}\nonumber\\
&= 2\int_{t- C \phi(r)}^{t-C\phi(r)+h} \int_{B_1} p_s^{B_1}(x,w) J(dw,B_2) \, 	ds \q \mbox{(by \eqref{e:levy-jump}).} \label{e:jl1}
\end{align}
By Proposition \ref{p:hke} and Remark \ref{r:conthk}(b), we have that
\be \label{e:jl2}
\wt u_h (t,x) :=2 \int_{t- C \phi(r)}^{t-C\phi(r)+h} \int_{B_1} p_s^{B_1}(x,w) J(dw,B_2) \, 	ds \q \mbox{is continuous in $(C\phi(r),\infty) \times B_1$.}
\ee

By \eqref{e:regphi}, we choose $A>1, \kappa \in (0,1)$ such that 
\be\label{e:jl3}
 \phi(Ar) > 2 \phi(r), \mbox{ and } \phi(\kappa r)< (C_2-C_1) \phi(r)/2 \q\mbox{for all $r > 0$.}
\ee
Let $\delta_N \in (0,1)$ denote the constant in \ref{NLE}. For any ball $B'=B(x',r') \subset B_1$ such that 
$B(x', A \delta_N^{-1} r') \subset B_1$, we have
\be \label{e:jl4}
p^{B_1}_s(x',w) \ge p_s^{B(x', A \delta_N^{-1} r')}(x',w) \gtrsim \frac{1}{V(x',r')}, \q \mbox{for all $s \in [\phi(r'),2\phi(r')]$ and $w \in B(x',r')$.}
\ee
We use  \ref{NLE}, \eqref{e:regphi} and \ref{VD} to obtain the above estimate. 
Set
\be \label{e:jl5}
\delta:= \min \left(c_0, \kappa , (A \delta_N^{-1}+1)^{-1}\right).
\ee 
The constant $\delta \in (0,1)$ is chosen so that for any ball $B'=B(x',r') \subset B(x_1,\delta r_1)$ with $r' \le \delta r_1$, we have
\be \label{e:jl6}
(C \phi(r)+\phi(r'),x') \in (C_1 \phi(r),C_2 \phi(r)) \times B(x_1, c_0r_1), \mbox{ and } B(x', A \delta_N^{-1} r') \subset B(x_1,r_1).
\ee
 
For any ball $B(x',r') \subset B(x_1,\delta r_1)$ with $r' \le \delta r_1$, we have
\begin{align} \label{e:jl7}
\esssup_{B(x_1, c_0r_1) \times (C_1 \phi(r),C_2 \phi(r))} u_{\phi(r')}(t,x) & \ge 
\wt u_{\phi(r')}\left(C \phi(r)+ \phi(r'), x'\right) \q \mbox{(by \eqref{e:jl1}, \eqref{e:jl2}, and \eqref{e:jl6})} \nonumber \\
&= 2 \int_{\phi(r')}^{2 \phi(r')} \int_{B_1} p^{B_1}_s(x,w)\, J(dw, B_2) \, ds \q \mbox{(by \eqref{e:jl2})} \nonumber \\
&\gtrsim  \phi(r')  \frac{J(B(x',r') \times B_2)}{\mu(B(x',r'))} \q \mbox{(by \eqref{e:jl4} and \eqref{e:jl6}).}
\end{align}

Set $C'=(C_3+C_4)/2$. For any $r' \le \delta r_1$, we have
 \begin{align} \label{e:jl8}
 \essinf_{B(x_1, c_0r_1) \times (C_1 \phi(r),C_2 \phi(r))} u_{\phi(r')}(t,x) & \le 
 \wt u_{\phi(r')}\left(C' \phi(r),x_1\right) \q \mbox{(by \eqref{e:jl1} and \eqref{e:jl2})} \nonumber \\
 &\le \int_{(C'-C)\phi(r)}^{(C'-C)\phi(r)+h} \int_{B_1} p_s(x,w)\, J(dw,B_2)\,ds  \q \mbox{($p^{B_1} \le p$)} \nonumber\\
& \lesssim \phi(r') \frac{J(B_1 \times B_2)}{\mu(B_1)} \q \mbox{(by \ref{UE}, \ref{VD}, and \eqref{e:regphi})}
 \end{align}
 The conclusion follows from \eqref{e:jl7}, \eqref{e:jl8} and \ref{PHI}.
 \qed

\noindent \emph{Proof of Theorem \ref{t:jkexists}.}
	Assume to the contrary that $J$ is not absolutely continuous with respect to $\mu \otimes \mu$ on  $(M\times M) \setminus d_M$, where $d_M= \set{(x,x): x \in M}$ denotes the diagonal in $M$.

	Let $\rho$ be the metric on $M \times M$ defined by $\rho((x_1,y_1),(x_2,y_2))= \max(d(x_1,x_2) , d(y_1,y_2))$. It is easy to verify that the product measure $\mu \otimes \mu$ satisfies the volume doubling property \eqref{VD}  on  the product space $(M \times M, \rho)$.
	For $(x_1,x_2)\in M \times M$, let $B_\rho((x_1,x_2),r)$ denote the open ball of radius $r$ in the metric $\rho$ centered at $(x_1,x_2)$.
	
	By the inner regularity of $J$, there exists $K$ a compact subset of $(M\times M) \setminus d_M$ such that $J(K)>0$ and $\mu \otimes \mu (K)=0$. Let $\delta>0$ be the constant in the statement of Lemma \ref{l:jumpav}. By the compactness of $K$, we can cover $K$ with finitely many sets of the form $\set{B(x,\delta d(x,y/4)) \times B(y,\delta d(x,y)/4): (x,y)\in K}$. Therefore, there exists $(x,y) \in K$ such that $\wt{K}=K \cap \left[B(x,\delta d(x,y/4)) \times B(y,\delta d(x,y)/4)\right]$ satisfies $J(\wt K)>0$ and $(\mu \otimes \mu) (\wt K)=0$.
%	Therefore, it is enough to assume that $K \subset B(x,\delta d(x,y/4)) \times B(y,\delta d(x,y)/4)$ for some $x,y \neq 0$. 

		By the regularity of $\mu \otimes \mu$, for any $\epsilon>0$, there exists an open set  $K_{\epsilon} \subset B(x,\delta d(x,y/4)) \times B(y,\delta d(x,y)/4) , \wt K \subset K_{\epsilon}$ such that $\mu \otimes \mu (K_{\epsilon})<\epsilon$. By the 5B-covering lemma \cite[Theorem 1.2]{Hei}, there exists balls $B_\rho((x_i,y_i), \rho_i) \subset K_{\epsilon}, i \in I$ such that $\rho_i \le \delta d(x_i,y_i)/4$ for all $i \in I$, $\bigcup_{i \in I} B_\rho((x_i,y_i), \rho_i) = K_{\epsilon}$ and $B_\rho((x_i,y_i),  \rho_i)/5), i \in I$ are pairwise disjoint. 
		Hence, we have
		\begin{align*}
		J(\wt K) & \le J(K_\eps) \le \sum_{i \in I} J(B_\rho((x_i,y_i), \rho_i) )= \sum_{i\in I} J\left(B((x_i,\rho_i)\times B(y_i,\rho_i))\right)  \\
		&\lesssim \sum_{i\in I}  \frac{ \mu(B((x_i,\rho_i))}{\mu(B(x, d(x,y)/4)} J\left(B(x,d(x,y)/4)\times B(y_i,\rho_i)\right)  \q \mbox{ (by Lemma \ref{l:jumpav})}\\
		&\lesssim  \sum_{i\in I}  \frac{ \mu(B((x_i,\rho_i))\mu(B((y_i,\rho_i))}{\mu(B(x,d(x,y)/4)\mu(B(y,d(x,y)/4))} J\left(B(x,d(x,y)/4)\times B(y,d(x,y)/4)\right) \\
		& \qq  \mbox{(by Lemma \ref{l:jumpav} and symmetry of $J$) }\\
		&\lesssim\frac{ J\left(B(x, d(x,y)/4)\times  B(y,d(x,y)/4)\right)}{\mu(B(x, d(x,y)/4)\mu(B(y,d(x,y)/4))} \sum_{i \in I} (\mu \otimes \mu) (B_\rho((x_i,y_i), \rho_i/5))  \q \mbox{(by \ref{VD})}\\
			&\lesssim \frac{ J\left(B(x, d(x,y)/4)\times  B(y, d(x,y)/4)\right)}{\mu(B(x,d(x,y)/4)\mu(B(y,d(x,y)/4)} (\mu \otimes \mu) (K_{\epsilon})  \\
			& \q  \mbox{(since $B_\rho((x_i,y_i),  \rho_i)/5), i \in I$ are pairwise disjoint and $\bigcup_{i \in I} B_\rho((x_i,y_i), \rho_i) = K_{\epsilon}$)}\\
	&\lesssim \epsilon \frac{ J\left(B(x,d(x,y)/4)\times  B(y,d(x,y)/4 )\right)}{\mu(B(x,d(x,y)/4)\mu(B(y,d(x,y)/4)} \q \mbox{(since $(\mu \otimes \mu) (K_{\epsilon}) < \eps$).}
		\end{align*}		
By letting $\epsilon \downarrow 0$, we obtain $J(\wt K)=0$, a contradiction.	
	\qed

\noindent  \\
Department of Mathematical Sciences, Tsinghua University, Beijing 100084, China.
\\and\\
Department of Mathematics, University of British Columbia,
Vancouver, BC V6T 1Z2, Canada. \\
liu-gh17@mails.tsinghua.edu.cn \sms\\

\noindent Department of Mathematics, University of British Columbia,
Vancouver, BC V6T 1Z2, Canada. \\
mathav@math.ubc.ca

\end{document}